\documentclass[11pt]{article}
\usepackage{amssymb}
\usepackage{amssymb}
\usepackage{amsmath}
\usepackage{latexsym}
\usepackage{theorem}
\usepackage{graphicx}
\usepackage{graphics}
\usepackage{epsfig}
 at 12 truept
\newtheorem{definition}{Definition}
\newtheorem{theorem}{Theorem}
\newtheorem{corollary}{Corollary}

\begin{document}
\date{}
\title{On the composition of integral operators acting in tempered Colombeau algebras}
\author{Alexei Filinkov\\
 School of Mathematical Sciences\\
University of Adelaide\\
Adelaide SA 5005, Australia\\
alexei.filinkov@adelaide.edu.au
\and Ian G. Fuss\\
School of Electrical and Electronic Engineering\\
University of Adelaide\\
Adelaide SA 5005, Australia\\
ian.fuss@adelaide.edu.au}
\maketitle
\date
\begin{abstract}
\noindent We show that the generalised composition of generalised integral operators is well defined on the space $\mathcal{G}_\tau$  Colombeau algebras of tempered generalised functions. 
\end{abstract}
\maketitle
\section{Introduction}
The extension of the classes of differential and integral equations that can be rigorously set and solved is seminal to current mathematics \cite{Grosser01} and vital to its application in diverse domains \cite{Ober16}. In this paper we continue the study of the field of generalised integral operators that was commenced in the context of Sobolev and Schwartz spaces of generalised functions \cite{Treves,Vladimirov}, ultradistributions \cite{Komatsu, Lozanov2007b} and continued recently within the spaces of Colombeau algebras of generalised functions \cite{Bernard2006a, Bernard2006b, Del2005, Del2007, Scarpalezos}. 

In part the motivation for the following analysis arises from a requirement in physics to be able to compose generalised integral operators\footnote{For a quick understanding of the need to generalise operator composition for quantum theory see section A.1 of Appendix A of Alexander Stottmeister's thesis \emph{On the Embedding of Quantum Field Theory on Curved Spacetimes into Loop Quantum Gravity} \cite{StottTh}} \cite{Folland,Stott16} and the fact that it is not possible to compose generalised integral operators that act within the space of Schwartz distributions \cite{Schwartz}.  It has been shown that such compositions exist in the Colombeau algebras of compactly supported generalized functions \cite{Del2007}.  
In this paper we demonstrate via an extension of the Schwartz kernel theorem to the space of bounded linear operators $\mathcal{L}(\mathcal{G}_\tau,\mathcal{G}_\tau)$ on tempered Colombeau algebras that compositions of generalised integral operators are well defined on the space $\mathcal{G}_\tau$.
More than this we use Hermite function expansions of ultradistributions to demonstrate that a countably infinite number of such compositions is well defined, hence we are able to show that compositions of exponentiated forms of these operators exist in the space $\mathcal{G}_{\tau_R}$ of Colombeau  tempered ultradistributions.

\section{Mapping from rapidly decreasing to tempered Colombeau algebras} 
It has been shown that operators $A$
\[
A:\mathcal{G}_S \rightarrow \mathcal{G}_\tau
\]
that are defined as
\begin{equation}\label{eqA}
A\phi = (A_\varepsilon)_\varepsilon (\phi_\varepsilon)_\varepsilon = \bigg(\int_{\mathbb{R}^n} K_\varepsilon(x,y)\phi_\varepsilon(y)dy \bigg)_\varepsilon \in \mathcal{G}_\tau
\end{equation}
where $\phi \in \mathcal{G}_S$ a Colombeau algebra of rapidly decreasing generalised functions, $\phi_\varepsilon \in S(\Omega)$ the Schwartz space of rapidly decreasing smooth functions,  and $\mathcal{G}_\tau$ is a Colombeau algebra of tempered generalised functions, are bounded linear operators: $A \in \mathcal{L}(\mathcal{G}_S, \mathcal{G}_\tau)$ \cite{Del2007}. Rigorous definitions of these spaces are given below.

\subsection{Simplified Colombeau algebras}

\subsubsection*{Colombeau algebra of rapidly decreasing generalised functions}
Let $\Omega$ be an open subset of $\mathbb{R}^d$ where $d\in\mathbb{N}_+$. Consider a smooth function $f\in C^\infty (\Omega)$ and denote
\[
\mu_{q,l}(f) :=\sup_{x\in\Omega,\, \vert\alpha\vert\leq l}\big(1 + \vert x\vert\big)^q\,\big\vert \partial^\alpha f(x)\big\vert\,,\quad\mathrm{where}\ q\in \mathbb{Z}\ \mathrm{and}\ l\in\mathbb{N}\,.
\]
The set
\begin{eqnarray*}
\lefteqn{
\mathcal{E}_\mathcal{S}(\Omega) := \Big\{ (f_\varepsilon)_\varepsilon\in \mathcal{S}(\Omega)^{(0,1]}\qquad \mathrm{such\ that}}\\ &&\qquad\qquad \forall q\,,l\in \mathbb{N}\  \exists n\in \mathbb{N}\,:\ \mu_{q,l}(f_\varepsilon) = O(\varepsilon^{-n})\ \mathrm{as}\ \varepsilon\to 0\Big\}
\end{eqnarray*}
is a sub-algebra of $\mathcal{S}(\Omega)^{(0,1]}$, where $\mathcal{S}(\Omega)$ is the Schwartz space of rapidly decreasing functions and
\[
\mathcal{S}(\Omega)^{(0,1]} = \big\{ (f_\varepsilon)_\varepsilon \ \mathrm{such\ that}\ \forall\varepsilon\in (0,1]\  f_\varepsilon\in \mathcal{S}(\Omega) \big\}
\]
is the set of nets. The set
\begin{eqnarray*}
\lefteqn{
\mathcal{N}_\mathcal{S}(\Omega) := \Big\{ (f_\varepsilon)_\varepsilon\in \mathcal{S}(\Omega)^{(0,1]}\qquad \mathrm{such\ that}}\\ && \qquad\qquad \forall q\,,l\in \mathbb{N}\  \forall p\in \mathbb{N}\,:\ \mu_{q,l}(f_\varepsilon) = O(\varepsilon^{p})\ \mathrm{as}\ \varepsilon\to 0\Big\}
\end{eqnarray*}
is an ideal in $\mathcal{E}_\mathcal{S}(\Omega)$. The factor-algebra
\[
\mathcal{G}_\mathcal{S}(\Omega) := \mathcal{E}_\mathcal{S}(\Omega)/\mathcal{N}_\mathcal{S}(\Omega)
\]
is referred to as the \emph{Colombeau algebra of rapidly decreasing generalised functions} (see \cite{Del2007} and references therein).

\subsubsection*{Colombeau algebras of tempered generalised functions}
Similarly, the factor-algebra
\[
\mathcal{G}_\tau(\Omega) := \mathcal{E}_\tau(\Omega)/\mathcal{N}_\tau(\Omega)
\]
is referred to as the \emph{Colombeau algebra of tempered generalised functions} \cite{Colomb84}. Here
\begin{eqnarray*}
\lefteqn{
\mathcal{E}_\tau(\Omega) := \Big\{ (f_\varepsilon)_\varepsilon\in \mathcal{O}_M(\Omega)^{(0,1]}\qquad \mathrm{such\ that}}\\ &&\qquad\qquad \forall l\in \mathbb{N}\  \exists q,\,n\in \mathbb{N}\,:\ \mu_{-q,l}(f_\varepsilon) = O(\varepsilon^{-n})\ \mathrm{as}\ \varepsilon\to 0\Big\}
\end{eqnarray*}
is a sub-algebra of $\mathcal{O}(\Omega)^{(0,1]}$, where
\[
\mathcal{O}_M(\Omega) := \Big\{ f\in C^\infty(\Omega)\ \mathrm{such\ that}\ \forall l\in \mathbb{N}\  \exists q\in \mathbb{N}\,:\ \mu_{-q,l}(f_\varepsilon) < \infty\Big\}
\]
is the algebra of smooth functions with slow growth (also known as the algebra of multiplicators); and
\begin{eqnarray*}
\lefteqn{
\mathcal{N}_\tau(\Omega) := \Big\{ (f_\varepsilon)_\varepsilon\in \mathcal{O}_M(\Omega)^{(0,1]}\qquad \mathrm{such\ that}}\\ &&\qquad\qquad \forall l\in \mathbb{N}\ \exists q\in \mathbb{N}\ \forall p\in \mathbb{N}\,:\ \mu_{-q,l}(f_\varepsilon) = O(\varepsilon^{p})\ \mathrm{as}\ \varepsilon\to 0\Big\}
\end{eqnarray*}
is an ideal in $\mathcal{E}_\tau(\Omega)$.
\subsubsection*{Colombeau algebras of tempered generalised ultradistributions}
Following \cite{Pilipovic2001} we define two types of Colombeau algebras of tempered generalised ultradistributions that correspond with the sets of generalised ultradistributions of Roumieu and Beurling type.
\begin{definition} 
The factor-algebra
\[
\mathcal{G}_{\tau_R}(\Omega) := \mathcal{E}_{\tau_R}(\Omega)/\mathcal{N}_{\tau_R}(\Omega)
\]
is defined as the Colombeau algebra of tempered generalised ultradistributions of Roumieu type.
Here 
\begin{eqnarray*}
\lefteqn{\mathcal{E}_{\tau_R}(\Omega) = \mathcal{E}_{\mathrm{exp}}^{\{M_p,N_p\}}(\Omega) :=\Big\{ (f_\varepsilon)_\varepsilon\in \mathcal{O}^{\{M_p\}}_{\mathrm{exp}}(\Omega)^{(0,1]}\qquad \mathrm{such\ that}}\\
&& \qquad\qquad \exists h,\,k >0\,:\ \nu_{h,M_p}(f_\varepsilon) = O(e^{N^\ast (k/\varepsilon)})\ \mathrm{as}\ \varepsilon\to 0\Big\} \,,
\end{eqnarray*}
is the set of generalised ultradistributions of Roumieu type, with 
\[
\mathcal{O}^{\{M_p\}}_{\mathrm{exp}}(\Omega) := \Big\{ f\in C^\infty(\Omega)\ \mathrm{such\ that}\ \exists h > 0\,:\ \nu_{h,M_p}(f_\varepsilon) < \infty\Big\}
\]
the Roumieu algebra of smooth functions of exponential growth; and
\begin{eqnarray*}
\lefteqn{\mathcal{N}_{\tau_R}(\Omega) = \mathcal{N}_{\mathrm{exp}}^{\{M_p,N_p\}}(\Omega) :=
\Big\{ (f_\varepsilon)_\varepsilon\in \mathcal{O}^{\{M_p\}}_{\mathrm{exp}}(\Omega)^{(0,1]}\qquad \mathrm{such\ that}}\\ && \qquad\qquad  \exists h>0,\,\forall k >0\,:\ \nu_{h,M_p}(f_\varepsilon) = O(e^{-N^\ast (k/\varepsilon)})\ \mathrm{as}\ \varepsilon\to 0\Big\}, 
\end{eqnarray*}
is an ideal in $\mathcal{E}_{\tau_R}(\Omega)$.
\end{definition}
Where as is customary in the theory of ultradistributions \cite{Komatsu}, we denote $M_p$ a sequence of positive numbers such that $M_0 = 1$ and
\begin{description}
\item[(M.1)] $M_p^2 \leq M_{p-1}M_{p-1}$ for any $p\in \mathbb{N}_+$;
\item[(M.2)] $M_p \leq c \,H^p\, M_q M_{p-q}$ for any $p\in \mathbb{N}_0\,, q\leq p$ and some $c, H\geq 1$;
\item[(M.3)] $\sum_{p=1}^{\infty} M_{p-1}/M_{p} < \infty$.
\end{description}
The sequence  $M_p^* := M_p/p!$ with $M_0^* = 1$,  the associated function
\[
M(\rho) := \sup_{p\in \mathbb{N}_0} \ln \frac{\rho^p}{M_p}\,,\quad \rho >0
\] 
and the growth function
\[
M^*(\rho) := \sup_{p\in \mathbb{N}_0} \ln \frac{\rho^p}{M^*_p}\,,\quad \rho >0\,.
\] 
For a smooth function $f\in C^\infty (\Omega)$ and we denote
\[
\nu_{h,M_p}(f) :=\sup_{x\in\Omega,\,\alpha,\beta\in\mathbb{N}_0^d } \frac{h^{\vert\alpha\vert + \vert\beta\vert} \big\vert x^\beta \partial^\alpha f(x)\big\vert}{M_{\vert\alpha\vert} \, M_{\vert\beta\vert}}\,,\quad\mathrm{where}\ h > 0\,.
\]

\begin{definition}
The factor-algebra
\[
\mathcal{G}_{\tau_B}(\Omega) := \mathcal{E}_{\tau_B}(\Omega)/\mathcal{N}_{\tau_B}(\Omega)
\]
is defined as the Colombeau algebra of tempered generalised ultradistributions of Beurling type.
Here 
\begin{eqnarray*}
\lefteqn{\mathcal{E}_{\tau_B}(\Omega) = \mathcal{E}_{\mathrm{exp}}^{(M_p,N_p)}(\Omega) := \Big\{ (f_\varepsilon)_\varepsilon\in \mathcal{O}^{(M_p)}_{\mathrm{exp}}(\Omega)^{(0,1]}\qquad \mathrm{such\ that}}\\ && \qquad\qquad  \forall h>0,\, \exists k >0\,:\ \nu_{h,M_p}(f_\varepsilon) = O(e^{N^\ast (k/\varepsilon)})\ \mathrm{as}\ \varepsilon\to 0\Big\},
\end{eqnarray*}
is the set of generalised ultradistributions of Beurling type, with
\[
\mathcal{O}^{(M_p)}_{\mathrm{exp}}(\Omega) := \Big\{ f\in C^\infty(\Omega)\ \mathrm{such\ that}\ \forall h > 0\,:\ \nu_{h,M_p}(f_\varepsilon) < \infty\Big\}.
\]
the Beurling algebra of smooth functions of exponential growth; and
\begin{eqnarray*}
\lefteqn{\mathcal{N}_{\tau_B}(\Omega) = \mathcal{N}_{\mathrm{exp}}^{(M_p,N_p)}(\Omega) := \Big\{ (f_\varepsilon)_\varepsilon\in \mathcal{O}^{(M_p)}_{\mathrm{exp}}(\Omega)^{(0,1]}\qquad \mathrm{such\ that}}\\ && \qquad\qquad  \forall h,\, k >0\,:\ \nu_{h,M_p}(f_\varepsilon) = O(e^{-N^\ast (k/\varepsilon)})\ \mathrm{as}\ \varepsilon\to 0\Big\} ,
\end{eqnarray*}
is an ideal in $\mathcal{E}_{\tau_B}(\Omega)$.
\end{definition}

\subsubsection*{Generalised constants}
We will also use the factor-ring of \emph{generalised constants}:
\[
\bar{\mathbb{K}} := \mathcal{E}_M(\mathbb{K})/\mathcal{N}(\mathbb{K})\,,
\]
for $\mathbb{K} = \mathbb{C}, \mathbb{R}\ \mathrm{or}\ \mathbb{R_+}$, where
\[
\mathcal{E}_M(\mathbb{K}) := \Big\{ (C_\varepsilon)_\varepsilon\in \mathbb{K}^{(0,1]}\ \mathrm{such\ that}\  \exists n\in \mathbb{N}\ :\ \vert C_\varepsilon\vert  = O(\varepsilon^{-n})\ \mathrm{as}\ \varepsilon\to 0\Big\}
\]
and
\[
\mathcal{N}(\mathbb{K}) := \Big\{ (C_\varepsilon)_\varepsilon\in
\mathbb{K}^{(0,1]}\ \mathrm{such\ that}\  \forall p\in \mathbb{N}\ :\ \vert C_\varepsilon\vert  = O(\varepsilon^{p})\ \mathrm{as}\ \varepsilon\to 0\Big\}
\]
\subsection{Inclusions}
Note that we have the following inclusions:
\[
\mathcal{G}_\tau(\Omega)\subset \mathcal{G}_{\tau_R}(\Omega)\subset \mathcal{G}_{\tau_B}(\Omega)\,.
\]

We also note the inclusions
\[
S_{\tau_R}^\prime(\Omega) \subset \mathcal{G}_{\tau_R}(\Omega)\quad\mathrm{and}\quad S_{\tau_B}^\prime(\Omega) \subset \mathcal{G}_{\tau_B}(\Omega) \,,
\]
where we denote $S_{\tau_R}^\prime$ the space of ultradistributions of Roumieu type and $S_{\tau_B}^\prime$ the space of ultradistributions of Beurling type.
Indeed, for $\varphi\in S_{\tau_R}(\Omega)$ and $f\in S^\prime_{\tau_R}(\Omega)$ we have \cite{Lozanov2007a,Vindas}
\[
\varphi = \sum_n a_n h_n\quad\mathrm{and}\quad f = \sum_n b_n h_n\,
\]
where $h_n$ are Hermite functions, which form an orthonormal basis of $L^2 (\mathbb{R}^d)$ \cite{Arfken} and Hermite coefficients $a_n$ and $b_n$ satisfy  estimates
\[
\vert a_n\vert \leq e^{-M(\sqrt{n}h)}\quad\mathrm{and}\quad \vert b_n\vert \leq e^{M(\sqrt{n}h)}\,.
\]
Define
\[
f_\varepsilon = \sum_n e^{-\varepsilon M^2(\sqrt{n}h)} b_n h_n \equiv \sum_n  f^\varepsilon_n h_n\,,
\]
where
\[
\vert f^\varepsilon_n\vert = \vert  e^{-\varepsilon M^2(\sqrt{n}h)} b_n\vert \leq C_\varepsilon\, e^{-M(\sqrt{n}h)}
\]
and therefore $f_\varepsilon\in S_{\tau_R}(\Omega)$ for each $\varepsilon > 0$. We observe that
\[
( f_\varepsilon )_\varepsilon = \Big( \sum_n f^\varepsilon_n h_n\Big)_\varepsilon \in \mathcal{G}_{\tau_R}(\Omega)\,,
\]
since
\[
\forall n,\ \varepsilon\qquad \vert f^\varepsilon_n\vert \leq  e^{M(\sqrt{n}h)}\,.
\]

\section{Mapping between tempered Colombeau algebras}

Since generalised integral operators of form (\ref{eqA}) 
\begin{equation*}
A\phi = (A_\varepsilon)_\varepsilon (\phi_\varepsilon)_\varepsilon = \bigg(\int_{\mathbb{R}^n} K_\varepsilon(x,y)\phi_\varepsilon(y)dy \bigg)_\varepsilon \in \mathcal{G}_\tau
\end{equation*}
are defined as bounded linear operators from $\mathcal{G}_S$ to $\mathcal{G}_\tau$, in order to compose such operators
we demonstrate that their extensions
\[
A:\mathcal{G}_\tau \rightarrow \mathcal{G}_\tau \, 
\]
can be well-defined. Such maps can be represented by nets $(A_\varepsilon )_\varepsilon$ of linear continuous maps 
\[
A = (A_\varepsilon)_\varepsilon \in \mathcal{L}(\mathcal{O}_M,\mathcal{O}_M)^{(0,1]} \, ,
\]
where these nets are defined to be of moderate growth if
\[
\forall \ell \in \mathbb{N} \ \exists \ (C_\varepsilon)_\varepsilon \in \mathcal{E}_M(\mathbb{R}_+), \ \exists \ p, q, \ell' \in \mathbb{N}
\]
such that
\[
\forall f \in \mathcal{O}_M \ \mu_{-p, \ell}(A_\varepsilon f) \leqslant C_\varepsilon \ \mu_{-q, \ell'}(f) \, .
\]
We then note that the net $(\phi_\varepsilon)_\varepsilon$ where $\phi_\varepsilon \in \mathcal{O}_M$ has an associated net
\[
\phi^\gamma_\varepsilon := \phi_\varepsilon e^{-\gamma |x|^2} \in S
\]
if
$\gamma \in \mathbb{R}_+$ with 
\[
\lim_{\gamma \rightarrow 0} \phi^\gamma_\varepsilon = \phi_\varepsilon \in S' \, .
\] 
In considering the nature of this limit it is helpful to note that $\mathcal{G}_\tau^\infty \cap S' = \mathcal{O}_M$ where $\mathcal{G}_\tau^\infty$ is the subspace of of regular elements of  $\mathcal{G}_\tau$ \cite{Del2008} and that the closure $\bar{S} = \mathcal{O}_M$ with convergence in $S'$.
We then define
\[
A\phi = (A_\varepsilon)_\varepsilon (\phi_\varepsilon^\gamma)_\varepsilon |_{\gamma = \varepsilon} = \bigg(\int_{\mathbb{R}^n} K_\varepsilon(x,y)\phi_\varepsilon^\gamma(y)dy \bigg)_\varepsilon\bigg|_{\gamma = \varepsilon} \in \mathcal{E}_\tau
\]
where we have used a double regularisation but for simplicity defined $\gamma = \varepsilon$. 
For any $\alpha$ there exists $q_1,q_2$ and $q_3 \in \mathbb{N}$ such that
\[
|\partial^\alpha \, A_\varepsilon \phi_\varepsilon^\gamma| \leqslant C(1 + |x|)^{q_1} \varepsilon^{-q_2} \gamma^{-q_3}|_{\gamma = \varepsilon} = C(1 +|x|)^{q_1}\varepsilon^{-q}
\]
where $q = q_2 +q_3$, therefore $ A_\varepsilon\in \mathcal{L}(\mathcal{O}_M,\mathcal{O}_M)$.

We note that for any $\alpha$ there exist $q_1, q_2\in \mathbb{N}$ such that
\begin{eqnarray*}
\lefteqn{\vert \partial_x^\alpha\, A_\varepsilon f_\varepsilon\vert = \big\vert \partial_x^\alpha\, A_\varepsilon (\phi_\varepsilon e^{-\varepsilon \vert y\vert^2})\big\vert}\\
&& \leq C\ \int_{\mathbb{R}^n} (1 + \vert x\vert)^{q_1}\, (1 + \vert y\vert)^{q_1} \vert f_\varepsilon (y)\vert \,dy\\
&& \leq C\ (1 + \vert x\vert)^{q_1}\,\int_{\mathbb{R}^n} (1 + \vert y\vert)^{q_1} \vert \phi_\varepsilon (y) e^{-\varepsilon \vert y\vert^2} (y)\vert \,dy \\
&& \leq C\ (1 + \vert x\vert)^{q_1}\, \sup_{y}\Big( (1 + \vert y\vert)^{-q_2}\vert \phi_\varepsilon (y)\vert \Big)\,\int_{\mathbb{R}^n} (1 + \vert y\vert)^{q_1 + q_2}  e^{-\varepsilon \vert y\vert^2} (y)\vert \,dy \\
&& \leq C\ (1 + \vert x\vert)^{q_1}\, \mu_{- q_2, 0}(\phi_\varepsilon)\,\varepsilon ^{-(q_1 + q_2)/2}\,,
\end{eqnarray*}
therefore for any $\alpha$
\[
(1 + \vert x\vert)^{-q_1}\,\big\vert \partial_x^\alpha\, A_\varepsilon (\phi_\varepsilon e^{-\varepsilon \vert y\vert^2})\big\vert\leq C\,\varepsilon ^{-q} \, \mu_{- q_2, 0}(\phi_\varepsilon)
\]
and thus for any $l$ there exist $p$ and $q^{\prime}$ such that
\[
\mu_{-p, l}\big( A_\varepsilon (\phi_\varepsilon e^{-\varepsilon \vert \cdot\vert^2})\big)\leq C\,\varepsilon ^{-q} \, \mu_{- q^{\prime}, 0}(\phi_\varepsilon)
\]
Now for any $\phi = (\phi_\varepsilon)_\varepsilon\in \mathcal{G}_\tau$ we define
\begin{eqnarray}\label{A}
A\phi &:=& \Big( A_\varepsilon (\phi_\varepsilon e^{-\varepsilon \vert \cdot\vert^2}) \Big)_\varepsilon + \mathcal{N}_\tau\\
&=& \bigg(\int_{\mathbb{R}^n}  K_\varepsilon (x,y) \phi_\varepsilon (y) e^{-\varepsilon \vert y\vert^2}dy \bigg)_\varepsilon + \mathcal{N}_\tau\nonumber
\end{eqnarray}
and we have that
\[
A \in \mathcal{L}(\mathcal{G}_\tau,\mathcal{G}_\tau) \, .
\]
\section{Composition of generalised integral operators on tempered Colombeau algebras} 
\begin{theorem}
Let generalised integral operators $A_1\,, A_2$ be defined by formula (\ref{A}), so that $A_1\,, A_2\in \mathcal{L}(\mathcal{G}_\tau,\mathcal{G}_\tau)$. Their composition $A_2 \circ A_1\in \mathcal{L}(\mathcal{G}_\tau,\mathcal{G}_\tau)$ is a generalised integral operator with the kernel
\[
K_\varepsilon (x,y) = \int_{\mathbb{R}^n} K^2_\varepsilon (x,z) K^1_\varepsilon (z,y) e^{-\varepsilon \vert z\vert^2}dz \in \mathcal{G}_\tau (\mathbb{R}^{2n})\,.
\]
\end{theorem}
{\bf{Proof}} 
For any $\alpha, \beta$ there exist $q_1, q_2\in \mathbb{N}$ such that
\begin{eqnarray}\label{est}
\vert \partial_x^\alpha\partial_y^\beta\, K_\varepsilon (x,y) \vert &=& \big\vert \partial_x^\alpha\partial_y^\beta\, \int_{\mathbb{R}^n} K^2_\varepsilon (x,z) K^1_\varepsilon (z,y) e^{-\varepsilon \vert z\vert^2}dz\big\vert\nonumber\\
&& \leq C\ (1 + \vert x\vert)^{q_1}\,(1 + \vert y\vert)^{q_2}\, \int_{\mathbb{R}^n} (1 + \vert z\vert)^{q_1 + q_2}  e^{-\varepsilon \vert z\vert^2} \,dz\nonumber \\
&& \leq C\ (1 + \vert x\vert)^{q_1}\,(1 + \vert y\vert)^{q_2}\,\varepsilon ^{-(q_1 + q_2)/2}\,,
\end{eqnarray}
therefore $K_\varepsilon\in \mathcal{G}_\tau (\mathbb{R}^{2n})$.

Furthermore
\begin{eqnarray*}
\big(A_2 \circ A_1\big)\phi 
&=& \int_{\mathbb{R}^n}  K_\varepsilon (x,y) \phi_\varepsilon (y) e^{-\varepsilon \vert y\vert^2}dy \\
&=& \int_{\mathbb{R}^n}\Bigg[\int_{\mathbb{R}^n} K^2_\varepsilon (x,z) K^1_\varepsilon (z,y) e^{-\varepsilon \vert z\vert^2}dz\Bigg] \phi_\varepsilon (y) e^{-\varepsilon \vert y\vert^2}dy \\
&=& \int_{\mathbb{R}^n} K^2_\varepsilon (x,z)\Bigg[\int_{\mathbb{R}^n}  K^1_\varepsilon (z,y)\phi_\varepsilon (y) e^{-\varepsilon \vert y\vert^2}dy\Bigg]  e^{-\varepsilon \vert z\vert^2}dz \\
&=&A_2\big(A_1\phi\big)\,.
\end{eqnarray*}

Estimate (\ref{est}) implies the following extension. 
\begin{corollary}
Let generalised integral operator $A\in \mathcal{L}(\mathcal{G}_\tau,\mathcal{G}_\tau)$ be defined by formula (\ref{A}). Then $A^k$ is well-defined in $\mathcal{L}(\mathcal{G}_\tau,\mathcal{G}_\tau)$ for  any $k$ and the operator 
\[
e^A := I + \sum_{k=1}^\infty A^k/ k!
\]
is well-defined in $\mathcal{L}(\mathcal{G}_{\tau_R},\mathcal{G}_{\tau_R})$. 
\end{corollary}


\begin{thebibliography}{99}
\bibitem{Arfken} G. Arfken, H. Weber, F.E. Harris, Mathematical Methods for Physicists, seventh ed., Academic Press, Cambridge, Massachusetts, 2001.

\bibitem{Bernard2006a} S. Bernard, J.-F. Colombeau, A. Delcroix, Generalized integral operators and applications. Math. Proc. Cambridge Philos. Soc. 141 (2006) 521-–546.

\bibitem{Bernard2006b} S. Bernard, J.-F. Colombeau, A. Delcroix, Composition and
exponential of compactly supported generalized integral operators, Integral Transforms Spec. Funct. 17 (2006) 93--99.

\bibitem{Colomb84} J.-F. Colombeau, Elementary Introduction to New Generalized Functions, North-Holland, Amsterdam, 1985.

\bibitem{Del2005} A. Delcroix, Generalized integral operators and Schwartz kernel type theorems, J. Math. Anal. Appl. 306 (2005) 481-–501.

\bibitem{Del2007} A. Delcroix, Kernel theorems in spaces of tempered generalized functions. Math. Proc. Cambridge Philos. Soc. 142 (2007) 557-–572.

\bibitem{Del2008} A. Delcroix,  A new approach to temperate generalized Colombeau functions. Publ. Inst. Math. (Beograd) (N.S.) 84(98) (2008) 109-–121.

\bibitem{Folland} G.B. Folland, Harmonic Analysis in Phase Space, Princeton University Press, Princeton, 1989.

\bibitem{Grosser01} M. Grosser, M.Kunzinger, M. Oberguggenberger, R. Steinbauer, Geometric Theory of Generalized Functions with Applications to General Relativity, Kluwer Academic Publishers, Dordrecht, 2001.

\bibitem{Komatsu} H. Komatsu, Ultradistributions, I: structure theorems and a characterization, J. Fac. Sci. Univ. Tokyo, Sect. IA Math. 20 (1973) 25–-105.

\bibitem{Lozanov2007a} Z. Lozanov-Crvenkovi$\mathrm{\acute{c}}$, D. Peri$\mathrm{\breve{s}}$i$\mathrm{\acute{c}}$, Hermite expansions of elements of Gelfand Shilov spaces in quasianalytic and non quasianalytic case. Novi Sad J. Math. 37 (2007) 129–-147.

\bibitem{Lozanov2007b}  Z. Lozanov-Crvenkovi$\mathrm{\acute{c}}$, D. Peri$\mathrm{\breve{s}}$i$\mathrm{\acute{c}}$,  Kernel theorems for the spaces of tempered ultradistributions. Integral Transforms Spec. Funct. 18 (2007) 699–-713.

\bibitem{Ober16}  H. Deguchi, M. Oberguggenberger, Propagation of singularities for generalized solutions to wave equations with discontinuous coefficients. SIAM J. Math. Analysis 48 (2016) 397--442.

\bibitem{Pilipovic2001} S. Pilipovic, D. Scarpalezos, Colombeau generalized ultradistributions, Math. Proc. Cambridge Philos. Soc., 130 (2001) 541--553.

\bibitem{Scarpalezos} D. Scarpalézos, Colombeau's generalized functions: topological structures; microlocal properties. A simplified point of view. II. Publ. Inst. Math. (Beograd) (N.S.) 76(90) (2004), 111-–125. 

\bibitem{StottTh} A. Stottmeister, On the embedding of quantum field theory on curved spacetimes into loop quantum gravity, 2015. Available from INIS: {http://inis.iaea.org/search/search.aspx?orig$\_$q=RN:47088807}

\bibitem{Stott16} A. Stottmeister, T. Thiemann, Coherent states, quantum gravity, and the Born-Oppenheimer approximation. III.: Applications to loop quantum gravity, J. Math. Phys. 57 (2016) 083509, https://doi.org/10.1063/1.4960823.
 
\bibitem{Schwartz} L. Schwartz, Sur l'impossibilite de la multiplications des distributions C. R. Acad. Sci. Paris, 239 (1954) 847--848.

\bibitem{Treves} F. Treves, Topological vector spaces, distributions and kernels. New York, Academic Press, 1967.

\bibitem{Vindas} D. Vu$\mathrm{\breve{c}}$kovi$\mathrm{\acute{c}}$,  J. Vindas,   Eigenfunction expansions of ultradifferentiable functions and ultradistributions in $\mathbb R^n$,
J. Pseudo-Differ. Oper. Appl. 7 (2016) 519--531.

\bibitem{Vladimirov} V.S. Vladimirov, Methods of the Theory of Generalized Functions, CRC Press, London, 2002.

\end{thebibliography}
\end{document}